\documentclass[12pt]{amsart}
\usepackage{amsmath,amsfonts,amssymb,amsrefs,fullpage,esint} 

\usepackage{amsthm}
\usepackage{tabularx}
\usepackage{tikz}
\usepackage{enumerate}
\newtheorem{theorem}{Theorem}

\DeclareFontFamily{U}{mathx}{}
\DeclareFontShape{U}{mathx}{m}{n}{<-> mathx10}{}
\DeclareSymbolFont{mathx}{U}{mathx}{m}{n}
\DeclareMathAccent{\widecheck}{0}{mathx}{"71}

\numberwithin{theorem}{section}

\newcommand{\R}{\mathbb{R}}
\newcommand{\T}{\mathbb{T}}
\newcommand{\N}{\mathbb{N}}
\newcommand{\Q}{\mathbb{Q}}
\newcommand{\Z}{\mathbb{Z}}
\newcommand{\W}{\mathbb{W}}
\renewcommand{\P}{\mathcal{P}}

\title{Strichartz inequalities: some recent developments}
\author{Jianhui Li}
\address{Department of Mathematics, Northwestern University, Evanston, IL-60208,
USA}
\email{jianhui.li@northwestern.edu}
\author{Zane Kun Li}
\address{Department of Mathematics, North Carolina State University, Raleigh, NC-27607,
USA}
\email{zkli@ncsu.edu}
\author{Po-Lam Yung}
\address{Mathematical Sciences Institute, Australian National University, Canberra, Australia \textit{and}  Department of Mathematics, The Chinese University of Hong Kong, Shatin, Hong
Kong }
\email{polam.yung@anu.edu.au  \phantom{aa} \textit{and} \phantom{aa} plyung@math.cuhk.edu.hk}

\begin{document}

\begin{abstract}
    Strichartz inequalities, originating from Fourier restriction theory, play a central role in the analysis of dispersive partial differential equations. They serve as a cornerstone for many subsequent developments. We survey some of them in memory of Strichartz, highlighting connections to recent developments in Fourier decoupling.
\end{abstract}

\dedicatory{In memory of Bob Strichartz}

\maketitle

\section{Introduction}

Strichartz inequalities have long been recognized as a standard way of capturing dispersion in spacetime. Their roots go all the way back to the Fourier restriction phenomena, first observed by E.M. Stein (unpublished), and has precedence in works of Fefferman \cite{MR257819}, Segal \cite{MR492892} and Tomas \cite{MR358216}. In his landmark paper \cite{MR512086}, Strichartz resolved fully the question of Fourier restriction to quadratic hypersurfaces in $\R^d$, of the form $S = \{x \in \R^d \colon R(x) = r\}$ where $R(x)$ is a quadratic polynomial of $d$ variables with real coefficients and $r$ is a constant. The cases where $S$ is a cone and $S$ is a paraboloid are of special interest; they correspond to dispersive estimates for the (acoustic) wave equation and the Schr\"{o}dinger equation respectively.

In this paper, we focus on the case of the Schr\"{o}dinger equation in $d + 1$ dimensional spacetime ($d \geq 1$) and the Strichartz inequality
\begin{equation} \label{eq:Str_paraboloid}
\|e^{i t \Delta} f(x)\|_{L^p(\R^{d+1})} \lesssim \|f\|_{L^2(\R^d)}, \quad p = \frac{2(d+2)}{d}
\end{equation}
first proved in \cite{MR512086}.
It has since been recognized that one can put a mixed norm on the left hand side of \eqref{eq:Str_paraboloid}:
\[
\|e^{i t \Delta} f(x)\|_{L^q_t(\R;L^p_x(\R^d))} \lesssim \|f\|_{L^2(\R^d)}, \quad \frac{2}{q}+\frac{d}{p} = \frac{d}{2}, \quad 2 \leq p,q \leq \infty, \quad (d,q,p) \ne (2,2,\infty).
\]
See Keel and Tao \cite{MR1646048} for the end-point estimate $(q,p) = (2,\frac{2 d}{d-2})$, Ginibre and Velo \cite{GV1992} and Yajima \cite{MR891945} for the non-endpoint cases, and Montgomery-Smith \cite{MR1600602} for showing the failure of the estimate when $(d, q, p) = (2, 2, \infty)$; the case $p = q$ corresponds to \eqref{eq:Str_paraboloid}. The above Strichartz inequality is the starting point of many important developments, such as Fourier extension inequalities for functions in $L^p$ ($p \ne 2$), bilinear and later multilinear improvements of \eqref{eq:Str_paraboloid} which went on to become indispensable tools in PDEs and harmonic analysis, attempts to find the best constants for \eqref{eq:Str_paraboloid}, and Strichartz inequalities for orthogonal systems. The relevant literature is huge and still growing. Instead of giving a detailed survey, we have opted to focus on two recent developments that we are more familiar with. We will discuss discrete and refined Strichartz estimates, from the point of view of Fourier decoupling.

Fourier decoupling inequalities provide ways to quantify interference among plane waves travelling in different directions. For our purpose, the relevant plane waves have frequencies lying in a small neighbourhood of the unit paraboloid $\mathfrak{P} := \{(\xi,|\xi|^2) \in \R^{d+1} \colon \xi \in [0,1]^d\}$. More precisely, suppose $N \gg 1$ and $P_{1/N}$ be a partition of $[0,1]^d$ into squares of side length $1/N$. For each $\theta \in P_{1/N}$, let $c_{\theta}$ be the center of $\theta$, $L_{\theta}(\xi) := |c_{\theta}|^2 + 2 c_{\theta} \cdot (\xi - c_{\theta})$ be the linearization of the map $\xi \mapsto |\xi|^2$ at $c_{\theta}$, and $\tau_{\theta}$ be the parallelepiped 
\begin{equation} \label{eq:tau_theta_def}
\tau_{\theta} := \{(\xi,\eta) \in \R^{d+1} \colon \xi \in 2 \theta,\, |\eta - L_{\theta}(\xi)| \leq N^{-2}\};
\end{equation}
here $2\theta$ is a square with the same center as $\theta$ but twice the side length.
The parallelepiped $\tau_{\theta}$ is essentially the smallest parallelepiped of dimensions $N^{-1} \times \dots \times N^{-1} \times N^{-2}$ containing a $N^{-2}$ neighbourhood of the piece of paraboloid $\mathfrak{P}$ sitting above $2 \theta$.  Let $\{f_{\theta}\}$ be a collection of Schwartz functions on $\R^{d+1}$ so that the Fourier transform $\widehat{f_{\theta}}$ is supported in $\tau_{\theta}$ for all $\theta$; we think of them as plane waves travelling in different directions because the parallelepiped $\{\tau_{\theta}\}$ are transverse thanks to the curvature of the paraboloid $\mathfrak{P}$. It is clear by Plancherel and the finite overlapping property of $\{\tau_{\theta}\}$ that 
\[
\Big\| \sum_{\theta \in P_{1/N}} f_{\theta} \Big\|_{L^2(\R^{d+1})} \lesssim \Big( \sum_{\theta \in P_{1/N}} \|f_{\theta}\|_{L^2(\R^{d+1})}^2 \Big)^{1/2}.
\]
In a breakthrough paper \cite{MR3374964}, Bourgain and Demeter showed that this is almost true if $L^2$ is replaced by $L^p$, $2 \leq p \leq \frac{2(d+2)}{d}$. More precisely, they showed that under the same assumptions,
\begin{equation} \label{eq:BD}
\Big\| \sum_{\theta \in P_{1/N}} f_{\theta} \Big\|_{L^p(\R^{d+1})} \lesssim_{\varepsilon} N^{\varepsilon} \Big( \sum_{\theta \in P_{1/N}} \|f_{\theta}\|_{L^p(\R^{d+1})}^2 \Big)^{1/2}, \quad p = \frac{2(d+2)}{d}
\end{equation}
for all $\varepsilon > 0$ (and this value of $p$ is sharp; the inequality with the $N^{\varepsilon}$ factor fails for any larger value of $p$). This and similar inequalities have had powerful applications in the study of harmonic analysis, partial differential equations, analytic number theory and incidence geometry (see e.g. Demeter \cite[Chapter 13]{MR3971577}). In dimension $d = 1$, Li \cite{MR4276295} sharpened the result \eqref{eq:BD} by showing that it continues to hold with $N^{\varepsilon}$ replaced by $\exp(C\frac{\log N}{\log \log N})$ for some constant $C > 0$. Guth, Maldague and Wang \cite{GMW} further improved this result and showed that when $d = 1$, \eqref{eq:BD} continues to hold if $N^{\varepsilon}$ is replaced by $(\log N)^c$ for some finite but unspecified constant $c$, i.e.
\begin{equation} \label{eq:GMW}
\Big\| \sum_{\theta \in P_{1/N}} f_{\theta} \Big\|_{L^6(\R^{2})} \lesssim (\log N)^c \Big( \sum_{\theta \in P_{1/N}} \|f_{\theta}\|_{L^6(\R^{2})}^2 \Big)^{1/2}.
\end{equation} 
A $p$-adic analogue of a variant of \eqref{eq:GMW} from \cite{GMW} but with an explicit $c$ was then obtained by the last two authors and Shaoming Guo in \cite{GLY}.
These are all closely related to the Strichartz estimate; in particular, they imply rather easily a Strichartz estimate on $\R^{d+1}$ with the corresponding (non-optimal) loss in the regularity of the initial data. They are also highly relevant for the discrete Strichartz estimate, which we discuss next.

In \cite{MR1209299}, Bourgain initiated the study of the Strichartz estimate \eqref{eq:Str_paraboloid} in the periodic setting, where one has less dispersion (because the waves keep coming back when the underlying physical domain is periodic). In particular, he proved that on $\T^2 := (\R / \mathbb{Z})^2$,
\begin{equation} \label{eq:discrete_Str}
\Big\|\sum_{n=1}^N a_n e^{2 \pi i (n x + n^2 t)} \Big\|_{L^6(\T^2)} \lesssim_{\varepsilon} N^{\varepsilon} \Big(\sum_{n=1}^N |a_n|^2 \Big)^{1/2}
\end{equation}
for every $\varepsilon > 0$ 
which is the appropriate periodic analogue of \eqref{eq:Str_paraboloid} for $d = 1$ (note that \eqref{eq:discrete_Str} cannot hold with $\varepsilon = 0$, even if all $a_n = 1$, as was shown in \cite[Remark 2]{MR1209299}). He also conjectured in higher dimensions ($d \geq 2$) that 
\[
\Big\|\sum_{|n| \leq N} a_n e^{2 \pi i (n \cdot x + |n|^2 t)} \Big\|_{L^p(\T^{d+1})} \lesssim_{\varepsilon} N^{\varepsilon} \Big(\sum_{|n| \leq N} |a_n|^2 \Big)^{1/2}, \quad p = \frac{2(d+2)}{d}
\]
for every $\varepsilon > 0$.
This was subsequently confirmed two decades later by himself and Demeter as an application of decoupling for the paraboloid \eqref{eq:BD} in $\R^{d+1}$. Bourgain actually had a slightly sharper bound in \cite{MR1209299} in dimension $d=1$, where $N^{\varepsilon}$ in \eqref{eq:discrete_Str} was replaced by $\exp(C \frac{\log N}{\log \log N})$ for some constant $C > 0$, and Li's result in \cite{MR4276295} matches that. The improved decoupling inequality \eqref{eq:GMW} of Guth, Maldague and Wang then gave an improvement of \eqref{eq:discrete_Str} where $N^{\varepsilon}$ is replaced by $(\log N)^c$ for some finite constant $c$, and the best result to date is the one in \cite{GLY} which says that $c$ can be taken to be $2+\varepsilon$ for any $\varepsilon > 0$.

On the other hand, very recently, in dimension $d = 2$, Herr and Kwak \cite{HerrKwak} has succeeded in proving the optimal discrete Strichartz inequality: by recasting the problem as one in incidence geometry and coming up with an ingenious counting argument, they were able to prove that for any finite subset $S \subset \mathbb{Z}^2$,
\begin{equation} \label{eq:HK}
\Big\| \sum_{n \in S} a_n e^{2\pi i (n \cdot x + |n|^2 t)} \Big\|_{L^4(\T^3)} \lesssim (\log \#S)^{1/4} \Big(\sum_{n \in S} |a_n|^2 \Big)^{1/2}.
\end{equation}
This is remarkable because prior to their work, even when $S = [-N,N]^2 \cap \mathbb{Z}^2$, it was not even known if one can have the estimate with a finite power of $\log N$ on the right hand side for this particular choice of $S$. Furthermore, the sharp dependence in $\log \#S$ in \eqref{eq:HK} is essential in establishing the small data global well-posedness of the cubic nonlinear Schr\"{o}dinger equation in $H^s(\T^2)$ for any $s > 0$. Since \cite{HerrKwak} only appeared after the completion of the current article, we will focus on the case $d = 1$ in what follows; we note that the Fourier decoupling approach described in this article, while interesting on its own, is unlikely to give anything as sharp as \eqref{eq:HK}.

Another way in which Fourier decoupling interacted with Strichartz inequalities arose in the work of Du, Guth and Li \cite{MR3702674}. There, a refined and bilinear version of \eqref{eq:Str_paraboloid} (with $d = 1$) was used to resolve Carleson's problem about pointwise a.e. convergence of the solution of the Schr\"{o}dinger equation in $\R^{2+1}$ to the initial data. They proved that such pointwise convergence holds when the initial data is in $H^s(\R^2)$ with $s > 1/3$. Higher dimensional refined Strichartz inequalities were established subsequently in Du, Guth, Li and Zhang \cite{MR3842310}, which led to partial progress to Carleson's problem in higher dimensions. This higher dimensional problem was finally resolved in full by Du and Zhang \cite{MR3961084} using a fractal $L^2$ refinement of previous ideas. Refined Strichartz inequalities also played an important role in the study of the Falconer distance problem, as in Guth, Iosevich, Ou and Wang \cite{MR4055179}, and in Du, Ou, Ren and Zhang \cite{DORZ2023}. 
They are typically proved (as in \cite{MR4055179}) using a refined version of the decoupling inequality \eqref{eq:BD}; see also Demeter \cite{arXiv:2002.09525} for perspectives, particularly about multilinear versions of refined Strichartz estimates, and generalizations of the linear one.

In this article, in order to get the main ideas across without getting bogged down by technical details, we have chosen to use a few heuristics that are white lies over $\R$. All these heuristics are actually true if one is willing to work $p$-adically over $\Q_p$; see e.g. \cite{Li_Qp} for details. We describe these heuristics in Section~\ref{sect:heuristics}. The discrete and refined Strichartz inequalities are then discussed in Sections~\ref{sect:discrete} and \ref{sect:refined} respectively.

\bigskip

\noindent {\it Acknowledgments.} J. Li is supported partially by NSF grant DMS-1653264 and Z. Li by NSF grants DMS-2037851 and DMS-2311174. Yung was among the many fortunate students who benefited from a legendary REU program that Bob Strichartz ran for 20 years at Cornell University. He is grateful for Bob's warm hospitality, kind encouragement, generous teaching and inspiring passion for mathematics. He is supported by a Future Fellowship FT20010039 from the Australian Research Council. We also thank the American Institute of Mathematics for their support through the Fourier restriction community, and Shaoming Guo for discussions from which we benefited a lot. 

\section{Heuristics / White lies} \label{sect:heuristics}

\subsection{The Fourier transform} \label{sect:FT}

The Fourier transform on $\R^d$ is defined by 
\[
\widehat{f}(\xi) = \int_{\R^d} f(x) e^{-2\pi i x \cdot \xi} dx.
\]
Similarly, the inverse Fourier transform is given by
\[
\widecheck{g}(x) = \int_{\R^d} g(\xi) e^{2\pi i x \cdot \xi} d\xi.
\]

One can use this formula to compute explicitly the inverse Fourier transform of the characteristic function $1_{[0,1]}$ on $\R$. On the other hand, for heuristics, it turns out to be more helpful to think of $\widecheck{1}_{[0,1]}$ as equal to $1_{[0,1]}$; this is certainly not quite true, but this provides a good intuition for many purposes (its analog on $\Q_p$ would be actually true). As a result, if $\theta$ is a parallelepiped in $\R^d$ of dimensions $r_1 \times \dots \times r_n$, then heuristically one should think of $\widecheck{1_{\theta}}$ as equal to $e^{2\pi i x \cdot \xi_0} |\theta^*|^{-1} 1_{\theta^*}$, where here $\xi_0$ is a point in $\theta$ and $\theta^* := \{x \in \R^d \colon |x \cdot (\xi-\xi_0)| \leq 1 \text{ for all $\xi \in \theta$}\}$ is a dual parallelepiped to $\theta$ passing through the origin with dimensions $r_1^{-1} \times \dots \times r_n^{-1}$.

Similarly, if $\theta$ is a cube in $[0,1]^d$ of side length $R^{-1/2}$, $\xi_0 \in \theta$, and $f = e^{2\pi i x \cdot \xi_0}1_{\theta^*}$, then heuristically $e^{i t \Delta} f(x) 1_{[0,R]}(t)$ can be thought of as $e^{2 \pi i (-2 \pi t |\xi_0|^2 + x \cdot \xi_0)} 1_T(x,t)$ %
where $T$ is a parallelepiped with base $\theta^* \times \{0\}$, height $R$ and long side pointing in the normal direction to the paraboloid $\mathfrak{P}$ at $(\xi_0,|\xi_0|^2)$. Here 
\[
e^{it\Delta} f(x) := \int_{\R^d} \widehat{f}(\xi) e^{-4\pi^2 i t |\xi|^2} e^{2\pi i x \cdot \xi} d\xi
\]
is the solution to the Schr\"{o}dinger equation $\frac{1}{i}\partial_t u = \Delta u$ with initial data $u(x,0) = f(x)$.
As a result, 
\[
\|e^{i t \Delta} f\|_{L^p([0,R]^{d+1})} \simeq |T|^{\frac{1}{p}} = R^{\frac{d+2}{2p}},
\]
which is comparable to $\|f\|_{L^2} = |\theta^*|^{1/2} = R^{d/4}$ if $p = \frac{2(d+2)}{d}$, the exponent in the Strichartz inequality \eqref{eq:Str_paraboloid}.

\subsection{The uncertainty principle}
We will use the following uncertainty principle, which again is only a heuristic
in $\R$ but is actually true in $\Q_{p}$. It asserts that if $\widehat{f}$
is supported in a spatial (or spacetime) parallelpiped $\tau$, then there exists
a tiling of space or spacetime (that is, $\R^d$ or $\R^{d + 1}$)
by translates of its dual parallelpiped $\tau^{\ast}$ such that $|f|$ is constant
on each such translate of $\tau^{\ast}$.

\subsection{Wave packet decomposition} \label{subsect:wp}

Suppose now $\delta \in \N^{-1}$. Recall $P_{\delta}$ is a partition of $[0,1]^d$ into a disjoint union of cubes of side length $\delta$. For $\theta \in P_{\delta}$, let $\T_{\theta}$ be a tiling of $\R^{d+1}$ by translates of $\tau_{\theta}^*$, which we think of as tubes. We write $\T_{\delta} := \bigsqcup_{\theta \in P_{\delta}} \T_{\theta}$. Every $T \in \T_{\delta}$ determines a unique $\theta$ such that $T \in \T_{\theta}$; we denote this $\theta$ by $\theta(T)$. A wave packet adapted to a tube $T \in \T_{\delta}$ is a function whose Fourier transform is supported on $2 \tau_{\theta(T)}$ and is supported on $T$. Technically this is impossible on $\R^{d+1}$ unless the wave packet is zero; the best one can do is to have the wave packet decaying rapidly off $T$. Hence in reality, one should thicken each $T$ by a factor that is a small power of $R$, but for the sake of presentation we will sweep this under the rug.

Suppose now $F$ is a Schwartz function supported on $\bigcup_{\theta \in P_{\delta}} \tau_{\theta}$. We take a partition of unity $\{\eta_{\tau_{\theta}}\}_{\theta \in P_{\delta}}$ with respect to the covering $\bigcup_{\theta \in P_{\delta}} \tau_{\theta}$. For $\theta \in P_{\delta}$, let
\[
\P_{\tau_{\theta}} F := F * \widecheck{\eta_{\tau_{\theta}}}
\]
so that $\widehat{\P_{\tau_{\theta}} F} = \widehat{F} \, \eta_{\tau_{\theta}}$.
For $T \in \T_{\delta}$, we let%
\[
\P_T F := 1_T \P_{\tau_{\theta(T)}} F.
\]
Following our earlier heuristics, we will pretend that $\P_T F$ is a wave packet adapted to $T$ (it is not quite true that $\P_T F$ has Fourier support in $2\theta(T)$; one should really use a cut-off that has compact Fourier support in place of $1_T$, but let's not worry about that here now). Since 
\[
F = \sum_{\theta \in P_{\delta}} \P_{\tau_{\theta}} F
\]
and
\[
\P_{\tau_{\theta}} F = \sum_{T \in \T_{\theta}} \P_T F \quad \text{for all $\theta \in P_{\delta}$},
\]
we have
\[
F = \sum_{T \in \T_{\delta}} \P_T F,
\]
a decomposition of $F$ into wave packets at scale $\delta$.

Similarly, we take a partition of unity $\{\psi_{\theta}\}_{\theta \in P_{\delta}}$ of the covering $\bigcup_{\theta \in P_{\delta}} 2\theta$. Given Schwartz $f$ on $\R^d$ with $\widehat{f}$ supported on $[0,1]^d$ and $\theta \in P_{\delta}$, we let%
\[
\P_{\theta} f := f * \widecheck{\psi_{\theta}}.
\]
Given $\nu \in \delta^{-1} \Z^d$, let 
\[
\P_{\theta,\nu} f :=  1_{q_{\nu}}\P_{\theta}f
\]
where $q_{\nu}$ is the cube in $\R^d$ of side length $\delta^{-1}$ centered at $\nu$. Then
\[
f = \sum_{(\theta,\nu) \in P_{\delta} \times \delta^{-1} \Z^d}\P_{\theta,\nu} f,
\]
and $\P_{\theta,\nu} f$ is a function which is supported in $q_{\nu}$ and whose Fourier support is (morally speaking) in $2\theta$. This decomposition is useful when we compute $e^{i t \Delta}f(x) 1_{[0,\delta^{-2}]}(t)$: in fact, it can be arranged so that to each $(\theta,\nu) \in P_{\delta} \times \delta^{-1} \Z^d$, there is a unique $T_{\theta,\nu} \in \T_{\theta}$ that lies in $\{t \geq 0\}$ and intersects $\{t = 0\}$ at $q_{\nu}$. Then $e^{it \Delta} (\P_{\theta,\nu} f)(x) 1_{[0,\delta^{-2}]}(t)$ is essentially a wave packet associated to the tube $T_{\theta,\nu}$.

\section{Strichartz inequality for the torus} \label{sect:discrete}

In this section we will discuss how one can bound
\begin{equation} \label{eq:exp_sum_L6T}
\Big\|\sum_{n=1}^N a_n e(n x + n^2 t) \Big\|_{L^6(\T^2)}
\end{equation}
by $\|a_n\|_{\ell^2} = \Big(\sum_{n=1}^N |a_n|^2 \Big)^{1/2}$ up to powers of $\log N$, where $e(t) := e^{2 \pi i t}$. This was a corollary of the improved decoupling inequality \eqref{eq:GMW} by Guth, Maldague and Wang \cite{GMW}; they derived \eqref{eq:GMW} by a powerful `high-low' method, that has since yielded success in a range of different problems in Fourier analysis. Below we illustrate the high-low method from \cite{GMW} by specializing the context of exponential sums, using the heuristics in Section~\ref{sect:heuristics}. We refer the readers to \cite{GM} where Guth and Maldague modified the high-low method from \cite{GMW} and used it to obtain an amplitude dependent wave envelope estimate, that yields not only standard $L^6$ decoupling inequalities for the parabola, but also their small cap variants. See also a more detailed exposition in Johnsrude \cite{Ben}.

The number of terms in the exponential sum in \eqref{eq:exp_sum_L6T} is $N$, which will be fixed throughout; it will be convenient to denote by
\[
Q := [0,N^2]^2
\]
a square of side length $N^2$. We will first make a change of variables and use periodicity to rewrite \eqref{eq:exp_sum_L6T} as 
\[
\Big( \fint_{Q} \Big|\sum_{n=1}^N a_n e(\gamma(\frac{n}{N}) \cdot x) \Big|^6 dx \Big)^{1/6}
\]
where $\gamma(s) := (s,s^2)$ and $\fint_Q$ denotes an averaged integral $\frac{1}{|Q|} \int_Q$; henceforth we write $x \in Q$ as opposed to what we previously denoted by $(x,t)$. We also assume $N \in 2^{\N}$. For $\delta \in 2^{-\N}$, we write $P_{\delta}$ for the partition of $[0,1]$ into a union of $\delta^{-1}$ many dyadic intervals of length $\delta$. 
\subsection{Initial attempt}
For $I \in P_{\delta}$, let 
\[
f_I(x) = \sum_{\substack{1 \leq n \leq N \\ n/N \in I}} a_n e( \gamma(\frac{n}{N}) \cdot x ),
\]
so that 
\[
f_{[0,1]}(x) = \sum_{I \in P_{\delta}} f_I(x).
\]
The classical Cordoba-Fefferman $L^4$ estimate (see for example \cite[Proposition 3.3]{MR3971577}) says essentially\footnote{This is a small lie, that one can make rigorous in various ways. One way is to incorporate weights on the right hand side, that is $\sim 1$ on $Q$ and decays rapidly off $Q$. In the following, we systematically ignore such technicalities, in favour of the exposition of main ideas.} that 
\[
\|f_{[0,1]}\|_{L^4_{\text{avg}}(B_{\delta^{-2}})} \leq C \Big\| \Big( \sum_{I \in P_{\delta}} |f_I|^2 \Big)^{1/2} \Big\|_{L^4_{\text{avg}}(B_{\delta^{-2}})} = C \Big\| \sum_{I \in P_{\delta}} |f_I|^2 \Big\|_{L^2_{\text{avg}}(B_{\delta^{-2}})}^{1/2} 
\]
for any square $B_{\delta^{-2}}$ of side length $\delta^{-2}$, where $\|f\|_{L^p_{\text{avg}(B)}} := (\fint_B |f|^p dx )^{1/p}$. If the following assumption holds true: 
\begin{equation} \tag{Assumption} \label{eq:miracle}
\text{$\{|f_I|^2\}_{I \in P_{\delta}}$ forms an orthogonal family on $L^2_{\text{avg}}(B_{\delta^{-2}})$,}
\end{equation}
then from $\|f_{[0,1]}\|_{L^6_{\text{avg}}(B_{\delta^{-2}})} \leq \|f_{[0,1]}\|_{L^{\infty}(B_{\delta^{-2}})}^{1/3} \|f_{[0,1]}\|_{L^4_{\text{avg}}(B_{\delta^{-2}})}^{2/3}$ we have
\begin{equation} \label{eq:strong_conclusion}
\|f_{[0,1]}\|_{L^6_{\text{avg}}(B_{\delta^{-2}})} \leq C \|f_{[0,1]}\|_{L^{\infty}(B_{\delta^{-2}})}^{1/3} \Big( \sum_{I \in P_{\delta}} \| f_I^2 \|_{L^2_{\text{avg}}(B_{\delta^{-2}})}^2 \Big)^{1/6}.
\end{equation}
This conclusion sounds plausible at the scale $\delta = 1/N$, since then $B_{\delta^{-2}} = Q$ and the right hand side is bounded by
\[
\Big( \sum_{n=1}^N |a_n| \Big)^{1/3} \Big( \sum_{n=1}^N |a_n|^4 \Big)^{1/6} = \|a_n\|_{\ell^1}^{1/3} \|a_n\|_{\ell^4}^{2/3},
\]
which dominates, by H\"{o}lder's inequality, $\|a_n\|_{\ell^2}$, the expected value of the left hand side $\|f_{[0,1]}\|_{L^6_{\text{avg}}(Q)}$ up to $N^{\varepsilon}$ losses. Unfortunately, the \eqref{eq:miracle} is certainly false at scale $\delta = 1/N$, since for every $I \in P_{1/N}$, the function $|f_I|^2$ is a non-negative constant function (just the modulus squared of one of the coefficients $a_n$), and constant functions cannot be orthogonal to each other on $Q$.

\subsection{Quantifying constructive interference}
The above strategy turns out to be useful in estimating $f_{[0,1]}$ where constructive interference happens. In this and the next subsection, consider the special situation where all $a_n = 1$. 
In this case, we would be estimating the $L^6_{\text{avg}}(Q)$ norm of 
\begin{equation} \label{eq:f_def}
f(x) := \sum_{n=1}^N e(\gamma(\frac{n}{N}) \cdot x).
\end{equation}
The $L^{\infty}$ norm of this function on $Q$ is of course equal to $N$, attained at $x = (jN,0)$ for $j = 0,1,\dots,N-1$. These are places where the constructive interference between the different $e(\gamma(\frac{n}{N}) \cdot x)$ are the strongest. Another way to measure constructive interference is to use a square function. For a scale $\delta \in [1/N, 1]$, we consider
\[
g(x) := \sum_{I \in P_{\delta}} |\sum_{n/N \in I} e(\gamma(\frac{n}{N}) \cdot x)|^2.
\]
The $L^{\infty}$ norm of $g(x)$ is again achieved at the points $x = (jN,0)$ for $j = 0,1,\dots,N-1$; this is because constructive interference occurs between all waves $e(\gamma(\frac{n}{N}) \cdot x)$ with $n/N$ belonging to the same interval $I \in P_{\delta}$. In fact such constructive interference gives $\|g\|_{L^{\infty}(Q)} = \delta^{-1} (N \delta)^2 = N^2 \delta$ which increases as the scale $\delta$ increases. 

With the above it might make sense to think of the set
\[
\Omega := \{x \in Q \colon |f(x)| \in [N/2,N]\},
\]
as where constructive interference is the strongest. For $x \in \Omega$, we have, by Cauchy-Schwarz, that 
\[
N/2 \leq  \sum_{I \in P_{\delta}} |\sum_{n/N \in I} e(\gamma(\frac{n}{N}) \cdot x)| \leq \delta^{-1/2} g(x)^{1/2}
\]
so
\[
g(x) \geq \frac{N^2 \delta}{4} \quad \text{for all $x \in \Omega$}.
\]
If we further have $\delta \geq 8/N$, then the above inequality implies $g(x) \geq 2N$, from which we have
\[
g(x) \leq 2 (g(x)-N) \quad \text{for all $x \in \Omega$}.
\]
This domination only holds on the set $\Omega$; the power of it lies with the cancellation that the right hand side exhibits on the square $Q$ (which is larger than $\Omega$). In fact, 
\[
g(x) - N = \sum_{I \in P_{\delta}} \Big( \Big|\sum_{n/N \in I} e(\gamma(\frac{n}{N}) \cdot x) \Big|^2 - N \delta \Big) = \sum_{I \in P_{\delta}} \sum_{\substack{n/N \in I \\ m/N \in I\\ n \ne m}} e\Big((\gamma(\frac{n}{N}) - \gamma(\frac{m}{N})) \cdot x\Big)
\]
is a sum of an orthogonal family on $Q$ indexed by $I \in P_{\delta}$; this plays the role of the desired but false \eqref{eq:miracle} in the derivation of \eqref{eq:strong_conclusion}. To see this, note that for $I, I' \in P_{\delta}$ with $I \ne I'$, we have
\[
\fint_{Q}  \sum_{\substack{n/N \in I \\ m/N \in I\\ n \ne m}} e\Big((\gamma(\frac{n}{N}) - \gamma(\frac{m}{N})) \cdot x\Big) \overline{ \sum_{\substack{n'/N \in I' \\ m'/N \in I' \\ n' \ne m'}}  e\Big((\gamma(\frac{n'}{N}) - \gamma(\frac{m'}{N})) \cdot x\Big) } dx = 0,
\]
since the left hand side is the number of solutions to 
\[
\begin{cases}
n-m = n'-m' \\
n^2-m^2=n'^2-m'^2
\end{cases}
\]
with $n/N, m/N \in I$, $n'/N,m'/N \in I'$, $n \ne m$ and $n' \ne m'$ (and no such solution exists because when $n-m = n'-m' \ne 0$, the second equation in the system gives $n+m = n'+m'$ as well, which together with the first equation gives $n=n'$, $m=m'$, contradicting $I \ne I'$). 
With this orthogonality in hand, we deduce, as a result, that
\[
\begin{split}
\int_{\Omega} g(x)^2 dx &\leq \int_{Q} 4 (g(x)-N)^2 dx \\
&= \int_{Q} 4 \sum_{I \in P_{\delta}} \Big| \sum_{\substack{n/N \in I \\ m/N \in I\\ n \ne m}} e\Big((\gamma(\frac{n}{N}) - \gamma(\frac{m}{N})) \cdot x\Big) \Big|^2 dx \\
&\leq \int_{Q} 4 \delta^{-1} (N \delta)^4 dx
\end{split}
\]
where we applied Cauchy-Schwarz pointwise in the last inequality. This shows
\[
\Big( \frac{1}{|Q|} \int_{\Omega} g(x)^2 dx \Big)^{1/6} \leq 2^{1/3} N^{2/3} \delta^{1/2}.
\]
But since $\|f\|_{L^{\infty}(Q)} = N$, we have
\[
\Big( \frac{1}{|Q|} \int_{\Omega} |f(x)|^6 dx \Big)^{1/6} \leq N^{1/3} \Big( \frac{1}{|Q|} \int_{\Omega} |f(x)|^4 dx \Big)^{1/6}.
\]
If we somehow have access to Cordoba-Fefferman's inequality at scale $\delta$ on this subdomain $\Omega$ of $Q$, namely
\begin{equation} \label{eq:CF-fake} \tag{CF small scale}
\int_{\Omega} |f(x)|^4 dx \leq C \int_{\Omega} g(x)^2 dx,
\end{equation}
then the above shows that
\[
\Big( \frac{1}{|Q|} \int_{\Omega} |f(x)|^6 dx \Big)^{1/6} \lesssim N^{1/3} N^{2/3} \delta^{1/2} = N \delta^{1/2}.
\]
If we can take $\delta$ to be say $\frac{(\log N)^{2c}}{N}$, then we get
\begin{equation} \label{eq:high_firsttry}
\Big( \frac{1}{|Q|} \int_{\Omega} |f(x)|^6 dx \Big)^{1/6} \lesssim (\log N)^c N^{1/2}
\end{equation}
which is the correct bound we expect\footnote{The set $\Omega$ contains balls of radius $O(1)$ around the $N$ points $(jN,0)$, $j=0,\dots,N-1$, where $|f| \simeq N$, so the integral on the left side is at least $(\frac{1}{N^4} N^{1+6})^{1/6} = N^{1/2}$.} on the set $\Omega$. The problem here is of course that inequality \eqref{eq:CF-fake} does not really hold with $\delta$ as small as $\frac{(\log N)^{2c}}{N}$. Indeed, $\int_{\Omega} |f|^4 \geq (N/2)^4 |\Omega|$ whereas $\int_{\Omega} g^2 \leq \|g\|_{L^{\infty}}^2 |\Omega| \leq (N^2\delta)^2 |\Omega|$, so if \eqref{eq:CF-fake} was true, then we must have $\delta \geq 1/4C$. But this is a relatively minor problem and can be circumvented by using bilinear restriction, in a similar way by which one would prove Cordoba-Fefferman.

To summarize, on the set $\Omega$ where we have a lot of constructive interference, the square function $g$ is close to its maximum, and hence much bigger than the average value of $g$ over $Q$ (namely $N$). One can then dominate $g$ by $g-N$ on the set $\Omega$, and then appeal to orthogonality on the bigger square $Q$. Such an argument has been used by Vinh \cite{MR2838005} in studying $r$-rich points in incidence geometry.

\subsection{Constructive interference at a hierarchy of scales}
In the previous subsection we assumed all $a_n = 1$, defined $f$ by \eqref{eq:f_def}, and explained a way to obtain a bound for $\Big( \frac{1}{|Q|} \int_{\substack{x \in Q \\ |f(x)| \simeq N}} |f(x)|^6 dx \Big)^{1/6}$ where $Q := [0,N^2]^2$. In this subsection, we continue to assume $a_n = 1$ for all $n$, and explain how one might show
\begin{equation} \label{eq:L6_alpha}
\Big( \frac{1}{|Q|} \int_{\substack{x \in Q \\ |f(x)| \simeq \alpha}} |f(x)|^6 dx \Big)^{1/6} \lesssim (\log N)^c N^{1/2}
\end{equation}
for other values of $\alpha \in (0,N/2)$. Since $\alpha < N/2$, on the set where $|f(x)| \simeq \alpha$ there will be some destructive interference. It will be helpful to quantify the amount of destructive interference by adopting a multi-scale point of view. We fix $\alpha$ in the following discussion; only $\sim \log N$ many dyadic values of $\alpha$ matters, so if we can show that \eqref{eq:L6_alpha} we have our desired bound for $(\fint_Q |f|^6)^{1/6}$.

Let $1/N = \delta_J < \delta_{J-1} < \dots < \delta_0 = 1$ be a range of scales. As a first attempt, we compute $f$ successively as follows. First, we let 
\[
f_{I_J} = e(\gamma(\frac{n}{N}) \cdot x) \quad \text{if $I_J \in P_{\delta_J}$ and $n/N \in I_J$}.
\]
Then we compute the sum
\[
f_{I_{J-1}} = \sum_{I_J \subset I_{J-1}} f_{I_J}, \quad \text{for each} \quad I_{J-1} \in P_{\delta_{J-1}}.
\]
Then we compute
\[
f_{I_{J-2}} = \sum_{I_{J-1} \subset I_{J-2}} f_{I_{J-1}}, \quad \text{for each} \quad I_{J-2} \in P_{\delta_{J-2}}.
\]
We continue computing $f_{I_{J-3}}$, \dots, $f_{I_1}$ until we reach $f_{I_0} = f$. We ask, at a given point $x$, whether constructive interference has occurred when we compute the sum $f_{I_j}$, for $j = J-1, J-2, \dots, 0$. To do so, we form the relevant square functions 
\[
g_j(x) = \sum_{I_j \in P_{\delta_j}} |f_{I_j}(x)|^2  1_Q(x), \quad j = J, \dots, 0.
\]
We then define, for a small constant $\varepsilon$ to be determined, a set
\[
\Omega_{J-1} := \{x \in Q \colon g_{J-1}(x) \geq (1+\varepsilon) g_J(x)\}.
\]
This is the set where constructive interference occurs when we computed $f_{I_{J-1}}$ for most $I_{J-1} \in P_{\delta_{J-1}}$, and is the analog of $\Omega$ in the last subsection. Next, we define
\[
\Omega_{J-2} := \{x \notin \Omega_{J-1} \colon g_{J-2}(x) \geq (1+\varepsilon) g_{J-1}(x)\},
\]
\[
\Omega_{J-3} := \{x \notin \Omega_{J-1} \cup \Omega_{J-2} \colon g_{J-3}(x) \geq (1+\varepsilon) g_{J-2}(x)\},
\]
etc, until we reach 
\[
\Omega_{0} := \{x \notin \Omega_{J-1} \cup \dots \cup \Omega_1 \colon g_{0}(x) \geq (1+\varepsilon) g_{1}(x)\}.
\]
In short, $\Omega_j$ is the set of $x$ where constructive interference occurred in computing $f_{I_j}$ for most $I_j \in P_{\delta_j}$, and destructive interference occurred when waves are superimposed at all frequency scales finer than $\delta_j$. We also define
\[
L := Q \setminus (\Omega_{J-1} \cup \dots \cup \Omega_0),
\]
which is the set where we always had destructive interference. The estimate for $f$ on $L$ is easy, because $|f|$ is small there: we have $|f(x)|^2 = g_0(x) \leq (1+\varepsilon) g_1(x) \leq \dots \leq (1+\varepsilon)^J g_J(x)$ for $x \in L$, so
\[
\Big( \frac{1}{|Q|} \int_L |f(x)|^6 dx \Big)^{1/6} \leq (1+\varepsilon)^{J/2} \Big( \fint_{Q} g_J(x)^3 dx \Big)^{1/6} = (1+\varepsilon)^{J/2} N^{1/2}
\]
which is a bound of the same quality as \eqref{eq:high_firsttry} if $\varepsilon$ is chosen sufficiently small so that $(1+\varepsilon)^{J/2} \lesssim (\log N)^c$. To estimate $f$ on $\Omega_j$, for $j = J,J-1,\dots,0$, first note that the Fourier transform of $g_j(x)$ is essentially supported in the union of a family of parallelepipeds of sizes $\delta_j \times \delta_j^2$ centered at the origin. If we compute the `low-frequency' part of this function, by localizing to the ball of radius $\delta_{j+1}$ centered at the origin, what we get is essentially $g_{j+1}(x)$ (this is a small lie, that can be rigorously justified if we work $p$-adically where we have perfect a uncertainty principle; see \cite[Lemma 7.1]{GLY}). Nonetheless we have 
\[
g_j * \eta_{j+1} \leq g_{j+1}
\]
where $\eta_{j+1}(x) := \delta_{j+1}^{2} \eta(\delta_{j+1} x)$ for a Schwartz function $\eta$ with compact Fourier support (see Low Lemma 3.24 of \cite{GMW}). Thus on $\Omega_j$, we have
\[
g_j(x) \geq (1+\varepsilon) g_j * \eta_{j+1}(x)
\]
which is equivalent to
\[
g_j(x) \leq \frac{1+\varepsilon}{\varepsilon} (g_j(x) - g_j * \eta_{j+1}(x)).
\]
Thus
\[
\begin{split}
\frac{1}{|Q|} \int_{\Omega_j} g_j(x)^2 dx 
&\leq \frac{1+\varepsilon}{\varepsilon} \fint_{Q} (g_j(x) - g_j * \eta_{j+1}(x))^2 dx \\
&\lesssim \frac{1+\varepsilon}{\varepsilon} \frac{\delta_j}{\delta_{j+1}} \fint_{Q} \sum_{I_j \in P_{\delta_j}} |f_{I_j}(x)|^4 dx
\end{split}
\]
by the almost orthogonality of $|f_{I_j}|^2 1_{Q} - |f_{I_j}|^2 1_{Q}*\eta_{j+1}$ as $I_j$ varies over $P_{\delta_j}$ (see High Lemma 3.25 of \cite{GMW}). If we have Cordoba-Fefferman at scale $\delta_j$, then mimicking what we did for the case $\alpha \simeq N$ above (where every $x$ is in $\Omega = \Omega_{J-1}$), we have
\begin{equation} \label{eq:high_j}
\begin{split}
\Big( \frac{1}{|Q|} \int_{\substack{|f(x)| \simeq \alpha \\ x \in \Omega_j}} |f(x)|^6 dx \Big)^{1/6} &\lesssim \alpha^{1/3} \Big( \frac{1}{|Q|} \int_{\substack{|f(x)| \simeq \alpha \\ x \in \Omega_j}} |f(x)|^4 dx \Big)^{1/6}  \\
&\lesssim \alpha^{1/3} \Big( \frac{1}{|Q|} \int_{\Omega_j} g_j(x)^2 dx \Big)^{1/6} \\
&\lesssim \alpha^{1/3} \Big( \frac{1+\varepsilon}{\varepsilon} \frac{\delta_j}{\delta_{j+1}} \Big)^{1/6} \Big( \fint_{Q} \sum_{I_j \in P_{\delta_j}} |f_{I_j}(x)|^4 dx \Big)^{1/6}.
\end{split}
\end{equation}
One would like to pull out 2 copies of $f_{I_j}$ from the integral and bound them by $\|f_{I_j}\|_{L^{\infty}}$. The problem is that if $j < J$ then we do not have much control of $\|f_{I_j}\|_{L^{\infty}}$ for each $I_j \in P_{\delta_j}$. Fortunately, if $x \in \Omega_j$ is such that $|f(x)| \simeq \alpha$, then those wave packets from $f_{I_j} 1_Q$ whose $L^{\infty}$ norm is much bigger than $\|g_J\|_{L^{\infty}}/\alpha$ contribute only very little to $|f(x)|$. This can be seen as follows: indeed, for such $x$ we have, from $x \notin \Omega_J \cup \dots \cup \Omega_{j+1}$, that
\[
g_{j+1}(x) \leq (1+\varepsilon) g_{j+2}(x) \leq \dots \leq (1+\varepsilon)^{J-j-1} g_J(x)
\]
so
\[
g_j(x) = \sum_{I_j \in P_{\delta_j}} |\sum_{I_{j+1} \subset I_j} f_{I_{j+1}}(x)|^2  1_{Q}(x) \leq \frac{\delta_j}{\delta_{j+1}} g_{j+1}(x) \leq \frac{\delta_j}{\delta_{j+1}} (1+\varepsilon)^{J-j-1} \|g_J\|_{L^{\infty}}.
\]
From $|f(x)| \simeq \alpha$ we also get,
\[
\alpha \simeq |f(x)| = |\sum_{I_j \in P_{\delta_j}} f_{I_j}(x)|,
\]
and if we focus on the contribution from those wave packets of $f_{I_j}$ whose $L^{\infty}$ norm is bigger than $C \|g_J\|_{L^{\infty}}/\alpha$ for some constant $C$, then their total contribution, by Cauchy-Schwarz, is at most
\[
\sum_{\substack{I_j \in P_{\delta_j} \\ |f_{I_j}(x)| > C \|g_J\|_{L^{\infty}}/\alpha}} |f_{I_j}(x)|
\leq \sum_{\substack{I_j \in P_{\delta_j} \\ |f_{I_j}(x)| > C \|g_J\|_{L^{\infty}}/\alpha}} \frac{|f_{I_j}(x)|^2}{ C \|g_J\|_{L^{\infty}}/\alpha} \leq \frac{g_j(x) }{C \|g_J\|_{L^{\infty}}/\alpha} = \frac{\alpha}{C} \frac{\delta_j}{\delta_{j+1}} (1+\varepsilon)^{J-j-1},
\]
which is much smaller than $|f(x)|$ if $C$ is much larger than $\frac{\delta_j}{\delta_{j+1}} (1+\varepsilon)^J$. This suggests that maybe we can assume all wave packets from $f_{I_j}$ are indeed bounded by $C \|g_J\|_{L^{\infty}}/\alpha$. If that's the case, then we may return to \eqref{eq:high_j}, and as long as $\varepsilon$ is not too small (say $\varepsilon \simeq (\log N)^{-c}$ and the scales are sufficiently close to one another (so that $\delta_j/\delta_{j+1} \lesssim (\log N)^c$), we have
\[
\begin{split}
\Big( \frac{1}{|Q|} \int_{\substack{|f(x)| \simeq \alpha \\ x \in \Omega_j}} |f(x)|^6 dx \Big)^{1/6} & \lesssim (\log N)^c \|g_J\|_{L^{\infty}}^{1/3} \Big( \fint_{Q} \sum_{I_j \in P_{\delta_j}} |f_{I_j}(x)|^2 dx \Big)^{1/6} \\
&\simeq (\log N)^c N^{1/3 + 1/6} = (\log N)^c N^{1/2},
\end{split}
\]
again a bound of the same strength as \eqref{eq:high_firsttry}. This leads one to prune the wave packets in the actual argument by removing any wave packets with excessive $L^{\infty}$ norms when we define $f_{J-1}$, $f_{J-2}$, \dots, $f_1$, and consider the corresponding square functions $g_{J-1}$, $g_{J-2}$, \dots, $g_1$, as well as the sets $\Omega_{J-1}$, $\Omega_{J-2}$, \dots, $\Omega_1$, after the pruning process. We describe the pruning process below (modulo difficulties caused by tail terms from the uncertainty principle).

\subsection{Pruning of wave packets}
Let's return now to the case of general coefficients $\{a_n\}_{1 \leq n \leq N}$, and fix $\alpha > 0$. Motivated by the above considerations, we introduce scales 
\[
1/N \simeq \delta_J < \dots < \delta_0 = 1
\]
where 
\[
\delta_j := (\log N)^{-c j}
\]
for all $j$, so that $\frac{\delta_j}{\delta_{j+1}} = (\log N)^c$. This dictates what $J$ is; we then have $J \simeq \frac{\log N}{\log \log N}$. The $\varepsilon$ in the definitions of the high and low sets should then be taken to be 
\begin{equation} \label{eq:eps_def}
\varepsilon := (\log N)^{-1},
\end{equation}
so that $(1+\varepsilon)^J \lesssim 1$. The constant $C$ used in the wave packet pruning should then be $(\log N)^{\tilde{c}}$ for some large enough $\tilde{c}$, so that $C$ is much larger than $\frac{\delta_j}{\delta_{j+1}} (1+\varepsilon)^J$. The height at which we should prune the wave packets will be denoted by 
\begin{equation} \label{eq:lambda_def}
\lambda := (\log N)^{\tilde{c}} \frac{\|g_J\|_{L^{\infty}(Q)}}{\alpha}.
\end{equation}

We will use the following notation throughout: for a general $F$ defined on $\R^2$, we let
\begin{equation} \label{eq:FI_def}
F_I := \P_{\tau_I} (F 1_Q), \quad I \in P_{\delta}
\end{equation}
and
\begin{equation} \label{eq:FT_def}
F_T := \P_T (F 1_Q), \quad T \in \T_{\delta},
\end{equation}
where $\P_{\tau_I}$ and $\P_T$ are projections onto frequency $\tau_I$ and wave packet $T$ defined in Section~\ref{subsect:wp}.

Let
\[
f := \sum_{n=1}^N a_n e(\gamma(\frac{n}{N}) \cdot x).
\]
Our goal is to prove $\|f\|_{L^6_{\text{avg}}(Q)} \lesssim (\log N)^c \|a_n\|_{\ell^2}$ where $Q = [0,N^2]^2$. By a standard bilinearization argument, one can replace the left hand side of this inequality by
\[
\Big( \fint_Q \max_{\substack{I_1, I_1' \in P_{\delta_1} \\ I_1 \ne I_1'}} |f_{I_1} f_{I_1'}|^3 dx \Big)^{1/6}.
\]
($f_{I_1}$ and $f_{I_1'}$ are defined using \eqref{eq:FI_def}.) Using a broad/narrow dichotomy (sometimes also called the Bourgain-Guth method, see \cite[Chapter 7]{MR3971577} or \cite[Section 6]{GLY}), one is reduced to showing that
\begin{equation} \label{eq:level_set_est}
\alpha^6 \frac{|U_{\alpha}(f)|}{|Q|} \lesssim [(\log N)^c \|a_n\|_{\ell^2}]^6
\end{equation}
where
\[
U_{\alpha}(f) := \Big\{x \in Q \colon \max_{\substack{I_1, I_1' \in P_{\delta_1} \\ I_1 \ne I_1'}} |f_{I_1} f_{I_1'}|^{1/2}(x) \simeq \alpha, \Big( \sum_{I_1 \in P_{\delta_1}} |f_{I_1}(x)|^6 \Big)^{1/6} \lesssim (\log N)^{c'} \alpha \Big\}
\]
for some $c'$.

To do so, we first form the square function at the finest frequency scale:
\[
g_J := \sum_{I_J \in P_{\delta_J}} |f_{I_J}|^2.
\]
Next we prune away the wave packets in each $f_{I_J}$ whose sup norm exceeds $\lambda$ (recall $\lambda$ from \eqref{eq:lambda_def}).
We let 
\[
f_J := \sum_{T \in \T_{\delta_J} \colon \|f_T\|_{L^{\infty}} \leq \lambda} f_T
\]
where we pruned away those wave packets $f_T$ from $f$ (defined via \eqref{eq:FT_def}) whose sup norm exceeds $\lambda$. The square function of the next coarser scale is then defined by
\[
g_{J-1} := \sum_{I_{J-1} \in P_{\delta_{J-1}}} |f_{J,I_{J-1}}|^2,
\]
(again, $f_{J,I_{J-1}}$ is defined from $f_J$ using the notation in \eqref{eq:FI_def}) and we will be interested in whether constructive interference has occurred when we computed $g_{J-1}$. Motivated by our choice of $\Omega_{J-1}$ in the last subsection and the choice of $\varepsilon$ in \eqref{eq:eps_def}, this leads us to define the set 
\[
\Omega_{J-1} := \{x \in Q \colon g_{J-1}(x) \geq (1+(\log N)^{-1}) g_J(x)\}.
\]
We then prune the wave packets in each $f_J$ whose sup norm exceeds $\lambda$:
\[
f_{J-1} := \sum_{T \in \T_{\delta_{J-1}} \colon \|f_{J,T}\|_{L^{\infty}} \leq \lambda } f_{J,T},
\]
form the square function of $f_{J-1}$ at the next coarser scale:
\[
g_{J-2} := \sum_{I_{J-2} \in P_{\delta_{J-2}}} |f_{J-1,I_{J-2}}|^2,
\]
let
\[
\Omega_{J-2} := \{x \in Q \setminus \Omega_{J-1} \colon g_{J-2}(x) \geq (1+(\log N)^{-1}) g_{J-1}(x)\}
\]
and keep going. In other words, at each step, we will throw away a small collection of wave packets whose sup norms are too big, and then form the square function at the next coarser scale. We decide whether constructive interference has happened, and $x$ will be in $\Omega_j$ if one has constructive interference in forming $g_j$ at $x$ but not at any finer frequency scales. Finally, we still define the low set by
\[
L := Q \setminus (\Omega_J \cup \dots \cup \Omega_1).
\]
Heuristically, each $\Omega_j$ is the disjoint union of squares of side length $\delta_j^{-1}$.
This is an instance of the uncertainty principle in the sense that if a function has Fourier transform supported in a ball of radius $R$, its modulus is ``essentially constant" on balls of radius $1/R$ (see for example this blog post of Tao \cite{TaoUP} for more discussion). Note that $\widehat{g_{j}}$ is supported in a ball of radius $\delta_j$ and hence $g_j$ is locally constant on balls of radius $\delta_{j}^{-1}$. Similarly $\widehat{g_{j + 1}}$ is supported in a ball of radius $\delta_{j + 1}$ and hence $g_{j + 1}$ is locally constant on balls of radius $\delta_{j + 1}^{-1} > \delta_{j}^{-1}$. Thus $\Omega_j$ can be decomposed into squares of side length $\delta_{j}^{-1}$ on which each of $g_{j}$ and $g_{j + 1}$ are constant. 

One can then verify, broadly following what we did in the last subsection, that 
\[
\Big( \frac{1}{|Q|} \int_{\Omega_j} \max_{\substack{I_1, I_1' \in P_{\delta_1} \\ I_1 \ne I_1'}} |f_{j,I_1} f_{j,I_1'}|^2 \Big)^{1/4} \lesssim \lambda^{1/2} (\log N)^c \|a_n\|_{\ell^2}
\]
for $j = J-1, \dots, 1$ and 
\[
\Big( \frac{1}{|Q|} \int_{L} \max_{\substack{I_1, I_1' \in P_{\delta_1} \\ I_1 \ne I_1'}} |f_{1,I_1} f_{1,I_1'}|^3 \Big)^{1/6} \lesssim (\log N)^c \|a_n\|_{\ell^2}.
\]
By approximating $|f_{I_1} f_{I_1'}|$ by $|f_{j,I_1} f_{j,I_1'}|$ on $U_{\alpha}(f) \cap \Omega_j$, one can also show that
\[
\alpha^6 \frac{|U_{\alpha} (f) \cap \Omega_j|}{|Q|} \leq \alpha^2 \frac{1}{|Q|} \int_{\Omega_k} \max_{\substack{I_1, I_1' \in P_{\delta_1} \\ I_1 \ne I_1'}} |f_{j,I_1} f_{j,I_1'}|^2 
\]
for $j= J-1, \dots, 1$, and similarly
\[
\alpha^6 \frac{|U_{\alpha} (f) \cap L|}{|Q|} \leq \frac{1}{|Q|} \int_{L} \max_{\substack{I_1, I_1' \in P_{\delta_1} \\ I_1 \ne I_1'}} |f_{1,I_1} f_{1,I_1'}|^6. 
\]
By combining the above four inequalities, and remembering that $\lambda = (\log N)^{\tilde{c}} \frac{\|g_J\|_{L^{\infty}}}{ \alpha} \lesssim (\log N)^{\tilde{c}} \frac{\|a_n\|_{\ell^2}}{\alpha}$ we finish the proof of \eqref{eq:level_set_est}.

\subsection{Remarks} 
The above high-low method can be seen as a systematic way of exploiting the $L^2$ orthogonality. With the notations in the above subsection, the keys are 
\begin{enumerate}[1)]
    \item the low lemma, a pointwise estimate originating from orthogonality that says 
\[
\Big|\sum_{I \in P_{r}} f_I \Big|^2 * \widecheck{1}_{|\xi| \leq r} = \sum_{I \in P_{r}} |f_I|^2 * \widecheck{1}_{|\xi| \leq r}
\]
for $0 < r \leq 1$; and
    \item the high lemma, an $L^2$ orthogonality estimate which says
\[
\int \Big| \sum_{I_{j-1} \in P_{\delta_{j-1}}} |f_{I_{j-1}}|^2*\widecheck{1}_{|\xi| \geq \delta_{j}} \Big|^2 \leq \frac{\delta_{j-1}}{\delta_{j}} \sum_{I_{j-1} \in P_{\delta_{j-1}}} \int \Big| |f_{I_{j-1}}|^2*\widecheck{1}_{|\xi| \geq \delta_{j}} \Big|^2.
\]
Here the scales are given by $1/N = \delta_J < \delta_{J-1} < \dots < \delta_0 = 1$.
\end{enumerate}
The two together facilitated the passage between square functions at different scales. In a more recent article \cite{GM}, Guth and Maldague presented a streamlined version of the high-low method, where they were also able to use the scales $w_j := N^{-2} \delta_j^{-1}$, so that the full range of scales read $1/N^2 = w_0 < \dots < w_J = 1/ N = \delta_J < \dots < \delta_0 = 1$; they then have 
\begin{enumerate}[3)]
    \item an additional variant of the high lemma that reads
\[
\int \Big| \sum_{I \in P_{1/N}} |f_I|^2*\widecheck{1}_{|\xi| \geq w_{j-1}} \Big|^2 = \sum_{I_{j-1} \in P_{\delta_{j-1}}} \int \Big| \sum_{I \in P_{1/N}, \, I \subset I_{j-1}} |f_I|^2 * \widecheck{1}_{|\xi| \geq w_{j-1}} \Big|^2.
\]
\end{enumerate}
The smaller scales $w_j$ are relevant because if $I_j \in P_{\delta_j}$, then the convex hull of $\bigcup_{I \in P_{1/N}, \, I \subset I_j} \tau_I^*$, which can also be identified with the dual to $\bigcap_{I \in P_{1/N}, \, I \subset I_j} (\tau_I-\tau_I)$, is essentially a rectangle of size $w_j^{-1} \times N^2$. Such rectangles arise in the wave envelope estimate they prove in \cite{GM}. In fact, they also used the local constancy of $\sum_{I \in P_{1/N}, \, I \subset I_j} |f_I|^2$ on translates of $\bigcap_{I \in P_{1/N}, \, I \subset I_j} \tau_I^*$, which can be identified with the dual to the convex hull of $\bigcup_{I \in P_{1/N}, \, I \subset I_j} (\tau_I-\tau_I)$, when they applied a bilinear restriction estimate. Combined with a different pruning of $f$ introduced in \cite{GM}, this yields a different looking high-low proof of decoupling for the parabola, that also yields small cap decoupling without too much additional difficulties.

\section{Refined Strichartz inequalities and applications} \label{sect:refined}

In the last section, we saw how one could improve the discrete Strichartz inequality using improvements of the decoupling inequality. In this section, we will see how a refined Strichartz inequality can be deduced from a refined decoupling inequality. 

Note that all discussions in this section, including the statements and proofs of the theorems, are subject to the heuristics laid out in Section~\ref{sect:heuristics}; they would have to be modified appropriately if they were to hold rigorously on $\R^{d+1}$.

First we state one version of a refined Strichartz inequality from \cite[Theorem 3.1]{MR3842310}.

\begin{theorem}[Refined Strichartz inequality, version 1] \label{thm:RSI} Suppose that $f$ has frequency support in $[0,1]^d$, and $Q_1,Q_2, ...$ are disjoint lattice $R^{1/2}$-cubes in $[0,R]^{d+1}$ so that 
$$
\|e^{it\Delta}f(x)\|_{L^p(Q_j)} \sim \lambda, \quad \forall \, j.
$$
Suppose that these cubes are arranged in horizontal slabs of the form $\R^{d} \times \{t_0, t_0+R^{1/2}\}$, and that each slab contains $\sim \sigma$ cubes $Q_j$. Let $Y=\cup_{j} Q_j$. Then for any $\varepsilon > 0$,
\begin{equation}\label{eq:thm:RSI}
    \|e^{it\Delta}f(x)\|_{L^p(Y)}\lesssim_\varepsilon R^{\varepsilon}\sigma^{-(\frac{1}{2}-\frac{1}{p})}\|f\|_{L^2}, \quad p=\frac{2(d+2)}{d}.
\end{equation}
\end{theorem}

To make sense of this, note that \eqref{eq:thm:RSI} represents a gain over \eqref{eq:Str_paraboloid} if $\sigma$ is large, since the power of $\sigma$ is negative. However, the left hand side only involves the $L^p$ norm of $e^{it \Delta}f$ over a subset $Y$ of $[0,R]^{d+1}$, and $e^{it\Delta} f$ is required to be equally spread out in the cubes that make up $Y$. Thus the gain in the refinement \eqref{eq:thm:RSI} of \eqref{eq:Str_paraboloid} comes from the assumption that the solution $e^{it\Delta}f(x)$ is spread out in spacetime. 
One example when this happens is when $e^{it \Delta}f(x) 1_{[0,R]^{d+1}}(x,t)$ is the sum of $R^{d/2}$ many wave packets in the same direction with the same amplitude: if $\theta \in P_{R^{-1/2}}$, $\xi_0 \in \theta$ and
\[
f(x) = \sum_{x_0} e^{2\pi i (x-x_0) \cdot \xi_0} 1_{\theta^*}(x-x_0)
\]
where the sum is over a maximal collection of $R^{1/2}$ separated points $\{x_0\}$ in $[0,R]^d$, then one has $\|e^{it\Delta} f\|_{L^p(Q)} \simeq R^{
\frac{d+1}{2p}}$ 
for all cubes $Q$ of side length $R^{1/2}$ inside $[0,R]^{d+1}$, so one can take $Y = [0,R]^{d+1}$, and one can check that for $p = \frac{2(d+2)}{d}$, one has $\frac{d+1}{p} = - \frac{d}{2}(\frac{1}{2}-\frac{1}{p}) + \frac{d}{2}$, so
\[
\|e^{it\Delta}f(x)\|_{L^p([0,R]^{d+1})} \simeq R^{\frac{d+1}{p}} = (R^{\frac{d}{2}})^{-(\frac{1}{2} - \frac{1}{p})} (R^d)^{\frac{1}{2}} = (R^{\frac{d}{2}})^{-(\frac{1}{2} - \frac{1}{p})} \|f\|_{L^2},
\]
illustrating the optimality of \eqref{eq:thm:RSI}.

A similar example happens when $f = 1_{[0,1]^d}$. Then stationary phase shows that $|e^{it\Delta} f(x)| \simeq R^{-\frac{d}{2}}$ on the set where $|x|+|t| \simeq R$, thus $\|e^{it\Delta} f(x)\|_{L^p(Q)} \sim R^{-\frac{d}{2}} R^{\frac{d+1}{2p}}$ whenever $Q \subset \R^{d+1}$ is a cube of side length $R^{1/2}$ at a distance $\simeq R$ to the origin. Thus we can take $Y$ to be the union of such cubes, and then there each slab has $\sim \sigma = R^{d/2}$ cubes $Q$. We then have again
\[
\|e^{it\Delta}f(x)\|_{L^p(Y)} \simeq R^{-\frac{d}{2}} R^{\frac{d+1}{p}} = (R^{\frac{d}{2}})^{-(\frac{1}{2} - \frac{1}{p})} = (R^{\frac{d}{2}})^{-(\frac{1}{2}-\frac{1}{p})} \|f\|_{L^2},
\]
matching the inequality \eqref{eq:thm:RSI}. Note that in this example, heuristically $e^{it\Delta} f(x)$ is the sum of $R^{d/2}$ many wave packets through the origin, each pointing in a different ($R^{-1/2}$ separated) direction, with amplitude $R^{-\frac{d}{2}}$.

In either examples above, there was no $R^{\varepsilon}$ loss, contrary to the statement of Theorem \ref{thm:RSI}. It is not known whether the $R^{\varepsilon}$ loss in Theorem \ref{thm:RSI} actually occurs; such a loss was absent in the Strichartz inequality \eqref{eq:Str_paraboloid}, but some sort of loss is necessary in many other situations (such as the decoupling inequality \eqref{eq:BD}, the Fourier restriction conjecture, and the Kakeya conjecture at the respective critical exponents).

We remark that the frequency localization of the initial data $f$ on $[0,1]^d$ in Theorem \ref{thm:RSI} is merely a technicality. In fact, if \eqref{eq:Str_paraboloid} holds for all $f$ frequency supported on $[0,1]^d$, 
an application of the Littlewood-Paley square function estimate and parabolic rescaling implies that \eqref{eq:Str_paraboloid} holds for all $f$.

Next, we state another version of refined Strichartz inequality, from \cite[Theorem 4.4]{MR4055179}. Let $\Lambda_{R^{1/2}} = R^{1/2} \Z^d \cap [0,R]^d$. Recall from  Section~\ref{subsect:wp} that if $\nu \in \Lambda_{R^{1/2}}$ we wrote $q_{\nu}$ for the cube in $\R^d$ centered at $\nu$ with side length $R^{1/2}$. Moreover, for every $\theta \in P_{R^{-1/2}}$ we let $\T_{\theta}$ for a partition of $\R^{d+1}$ into translates of $\tau_{\theta}^*$; it was arranged so that for $(\theta,\nu) \in P_{R^{-1/2}} \times \Lambda_{R^{1/2}}$ there exists a unique tube $T_{\theta,\nu} \in \T_{\theta}$ that lies in $\{t \geq 0\}$ and intersects $\{t = 0\}$ at $q_{\nu}$.

\begin{theorem}[Refined Strichartz inequality, version 2] \label{thm:RSI2}
Suppose that $R \gg 1$, $\W \subset P_{R^{-1/2}} \times \Lambda_{R^{1/2}}$ and $f = \sum_{(\theta,\nu) \in \W} f_{\theta,\nu}$ is the sum of $W$ non-zero functions, all with comparable $L^2$ norms, so that each $f_{\theta,\nu}$ has Fourier support in $2\theta$ and is supported on $q_{\nu}$. Let $Y \subset [0,R]^{d+1}$ be a union of cubes $Q_j$ of side length $R^{1/2}$ in $\R^{d+1}$ where each $Q_j$ is contained in $\sim M$ tubes $T_{\theta,\nu}$ with $(\theta,\nu) \in \W$. Then
\[
    \|e^{it\Delta}f(x)\|_{L^p(Y)} \leq C_\varepsilon R^{\varepsilon} \Big( \frac{M}{W} \Big)^{\frac{1}{2}-\frac{1}{p}} \|f\|_{L^2}, \quad p = \frac{2(d+2)}{d}.
\]
\end{theorem}

Note that $M/W$ is always smaller than 1. It is the fraction of all wave packets from $e^{i t \Delta} f(x) 1_{[0,R]}(t)$ that pass through a cube $Q_j \subset Y$. Thus apart from the $R^{\varepsilon}$ loss, this strengthens the Strichartz inequality when $M/W$ is very small, i.e. when the wave packets of $e^{i t \Delta} f(x) 1_{[0,R]}(t)$ are very spread out.

It is folklore that Theorem~\ref{thm:RSI} follows from Theorem~\ref{thm:RSI2}; we thank Hong Wang for communicating this to the first author (see \cite[Section 2.3]{Li_thesis}). The proof of this implication is reproduced below. Its gist can also be found in Demeter \cite{arXiv:2002.09525}.

\begin{proof}[Proof of Theorem~\ref{thm:RSI}] %
The configurations of cubes in $Y$ in Theorems ~\ref{thm:RSI} and ~\ref{thm:RSI2} are different. We shall start by finding a significant subset of cubes in $Y$ that satisfies the assumptions in Theorem ~\ref{thm:RSI2}. The process is known as dyadic pigenholing. We sort tubes $T_{\theta,\nu} \in \T$ and cubes $Q_j \subset Y$ as follows.
\begin{itemize}
    \item For a dyadic number $\mu$, let $\W_\mu$ be the collection of pairs $(\theta,\nu)$ such that $\|f_{\theta,\nu}\|_{L^2} \sim \mu \|f\|_{L^2}$. By Strichartz inequality, the contributions from tubes in $\W_{\mu}$ with $\mu < R^{-C}$ are negligible. Moreover, $\|f_{\theta,\nu}\|_{L^2} \leq \|f\|_{L^2}$. $\W_{\mu}$ is empty when $\mu >1$. Therefore, We are only interested in sets $\W_\mu$ with $R^{-C} \leq \mu \leq 1$. 
    \item For each cube $Q_j \subset Y$, we have
    $$
    \|e^{it\Delta} f\|_{L^p(Q_j)} \leq \sum_{\mu} \Big \| \sum_{(\theta,\nu)\in \mathbb W_\mu }e^{it \Delta} f_{\theta,\nu} \Big\|_{L^p(Q_j)}.
    $$
    Since there are only $\sim \log R$ many interesting dyadic numbers $\mu$, the pigeonhole principle asserts that there exists a dyadic number $\mu(Q_j)$ such that $R^{-C} \leq \mu(Q_j) \leq 1$ and
    \begin{equation}\label{eq:dya_pigen1}
         \|e^{it\Delta} f\|_{L^p(Q_j)} \lesssim (\log R) \Big \| \sum_{(\theta,\nu)\in \mathbb W_{\mu(Q_j)} }e^{it \Delta} f_{\theta,\nu} \Big\|_{L^p(Q_j)}.
    \end{equation}
    Moreover, let $\mu$ be the most popular value of $\mu(Q_j)$ as $Q_j$ varies inside $Y$. Then $\gtrsim \frac{1}{\log R}$ portion of $Q_j$ in $Y$ satisfies $\mu(Q_j) = \mu$. We fix this value $\mu$ in the sequel and write $\mathbb{W} = \mathbb{W}_\mu$.
    \item We further sort cubes $Q_j \subset Y$ satisfying $\mu(Q_j) = \mu$ into a disjoint union $\bigsqcup_{M \in 2^{\N}} Y_M$ as follows. For each dyadic number $M$, let $Y_M$ be the set of all $Q_j$ with $\mu(Q_j) = \mu$ so that $Q_j$ intersects $\sim M$ many wavepackets $T_{\theta,\nu}$ with $(\theta,\nu) \in \W$. Since $Y_M$ is empty unless $1 \leq M \leq R^{d/2}$, there exists a dyadic number $M$ such that the collection $Y_{M}$ contains $\gtrsim \frac{1}{(\log R)^2}$ portion of $Q_j$ in $Y$. We fix this value $M$ in the sequel and write $Y'$ for the union of all $Q_j$ in $Y_M$.
\end{itemize}

Recall that $\|e^{it\Delta} f\|_{L^p(Q_j)} \sim \lambda$. Hence
\[
\|e^{it\Delta} f\|_{L^p(Y)}^p \lesssim \frac{|Y|}{|Q_j|}\lambda^p \lesssim (\log R)^2 \frac{|Y'|}{|Q_j|} \lambda^p \leq  (\log R)^{2} \sum_{Q_j \subset Y'} \|e^{it\Delta} f\|_{L^p(Q_j)}^p.
\]
In summary, we obtain
$$
 \|e^{it\Delta} f\|_{L^p(Y)} \lesssim (\log R)^{\frac{2}{p}} \Big(\sum_{Q_j \subset Y'} \|e^{it\Delta} f\|_{L^p(Q_j)}^p \Big)^{1/p} \lesssim (\log R)^{1+\frac{2}{p}} \Big\|\sum_{(\theta,\nu)\in \mathbb W }e^{it \Delta} f_{\theta,\nu}\Big\|_{L^p(Y')},
$$
where each $Q_j$ in $Y'$ intersects $\sim M$ many wavepackets $T_{\theta,\nu}$ with $(\theta,\nu) \in \mathbb W$, $\|f_{\theta,\nu}\|_{L^2}$ is dyadically a constant for $(\theta,\nu) \in \mathbb W$, and $|Y| \lesssim (\log R)^2 |Y'|$.

Let $W=|\mathbb{W}|$. By Theorem~\ref{thm:RSI2}, we obtain
\[
    \Big\|\sum_{(\theta,\nu)\in \mathbb W }e^{it \Delta} f_{\theta,\nu}\Big\|_{L^p(Y')} \lesssim_\varepsilon R^{\varepsilon} (\frac{M}{W})^{1/2 - 1/p } \|f\|_{L^2}.
\]
The proof of Theorem~\ref{thm:RSI} will be complete if we show that 
$$
\frac{M}{W} \lesssim (\log R)^C \sigma^{-1}.
$$

To prove the inequality, we count the number of $Q_j$ in $Y'$, denoted by $N$. Suppose that the number of horizontal slabs is $\rho$. By assumption,  the number of $Q_j \subset Y$ is $\sim \sigma \rho$. Using $|Y| \lesssim (\log R)^2 |Y'|$,
we have $\sigma \rho \lesssim (\log R)^2 N$.
On the other hand, each $T_{\theta,\nu}$ can pass through at most $\sim \rho$ $Q_j$'s. The number of interactions between $T_{\theta, \nu}$ and $Q_j$ is bounded above by $\sim \rho W$. On the other hand, the number of tubes $T_{\theta, \nu}$ intersecting $Q_j$ is $\sim$ $M$, and each $Q_j$ is disjoint. 
So $N \lesssim \frac{\rho W}{M}$. Putting together, we have $\sigma \rho \lesssim (\log R)^2 N \lesssim (\log R)^2 \frac{\rho W}{M}$,
which gives $\sigma \lesssim (\log R)^2 \frac{W}{M}$
 as desired.
\end{proof}

Next we need to prove Theorem \ref{thm:RSI2}. It relies on a refined version of the decoupling inequality \eqref{eq:BD}. 

Given $R \gg 1$, let $\T = \bigsqcup_{\theta \in P_{R^{-1/2}}} \T_{\theta}$. For $\W \subset \T$ we write $\W_{\theta} := \W \cap \T_{\theta}$. Our formulation of the refined decoupling inequality is inspired by the one in \cite{arXiv:2002.09525}.

\begin{theorem}[Refined decoupling inequality] \label{thm:refined_dec}
Let $\W \subset \T$ and $F = \sum_{T \in \W} F_T$ be a sum of wave packets from $\W$. Suppose $1 \leq M \leq R^{d/2}$ and $Y_M \subset \R^{d+1}$ be the union of all lattice cubes $Q$ of side length $R^{1/2}$ such that $\#\{T \in \W \colon Q \subset T\} \sim M$. Then for any $\varepsilon > 0$,
\begin{equation} \label{eq:refined_dec}
\|F\|_{L^p(Y_M)} \leq C_{\varepsilon} R^{\varepsilon} M^{\frac{1}{2}-\frac{1}{p}} \Big( \sum_{\theta \in P_{R^{-1/2}}} \|F_{\theta}\|_{L^p(\R^{d+1})}^p \Big)^{1/p}, \quad p = \frac{2(d+2)}{d}
\end{equation}
where $F_{\theta} := \sum_{T \in \W_{\theta}} F_T$.
\end{theorem}

When $M \sim R^{d/2}$, \eqref{eq:refined_dec} follows from \eqref{eq:BD} by using H\"{o}lder's inequality to estimate the $\ell^2$ norm over $\theta$ with an $\ell^p$ norm. Thus the main thrust of this inequality is the gain we obtain when we compute the $L^p$ norm of $F$ on a set that is incident to very few wave packets.

Below we first use Theorem~\ref{thm:refined_dec} to prove Theorem~\ref{thm:RSI2}, and then discuss a proof of Theorem~\ref{thm:refined_dec}.

\begin{proof}[Proof of Theorem~\ref{thm:RSI2}]
Let
\[
F(x) := e^{it \Delta} f(x) 1_{[0,R]}(t) = \sum_{(\theta,\nu) \in \W} e^{it \Delta} f_{\theta,\nu}(x) 1_{[0,R]}(t)
\]
which we think heuristically as the sum of wave packets. We are then in position to apply Theorem~\ref{thm:refined_dec} to estimate $\|F\|_{L^p(Y)}$: we should then let, for $\theta \in P_{R^{-1/2}}$,
\[
F_{\theta} := \sum_{\substack{\nu \in \Lambda_{R^{1/2}} \\ (\theta,\nu) \in \W}} e^{it \Delta} f_{\theta,\nu}(x) 1_{[0,R]}(t).
\]
In fact, if we let $
f_{\theta} := \sum_{\substack{\nu \in \Lambda_{R^{1/2}} \\ (\theta,\nu) \in \W}} f_{\theta,\nu},$
then
\[
\|F_{\theta}\|_{L^p(\R^{d+1})} = \|e^{it \Delta} f_{\theta}\|_{L^p(\R^d \times [0,R])} = \Big(\sum_{\substack{\nu \in \Lambda_{R^{1/2}} \\ (\theta,\nu) \in \W}} \|e^{it \Delta} f_{\theta}\|_{L^p(T_{\theta,\nu})}^p \Big)^{1/p} \lesssim \Big( \sum_{\substack{\nu \in \Lambda_{R^{1/2}} \\ (\theta,\nu) \in \W}} \|f_{\theta,\nu}\|_{L^2(\R^d)}^p \Big)^{1/p}
\]
so using the comparability between the non-zero $\|f_{\theta,\nu}\|_{L^2}$ and orthogonality, we have
\[
\begin{split}
\Big( \sum_{\theta \in P_{R^{-1/2}}} \|F_{\theta}\|_{L^p(\R^{d+1})}^p \Big)^{1/p} &\lesssim \Big( \sum_{(\theta,\nu) \in \W} \|f_{\theta,\nu}\|_{L^2(\R^d)}^p \Big)^{1/p} \\
&= W^{-(\frac{1}{2}-\frac{1}{p})} \Big( \sum_{(\theta,\nu) \in \W} \|f_{\theta,\nu}\|_{L^2(\R^d)}^2 \Big)^{1/2} = W^{-(\frac{1}{2}-\frac{1}{p})} \|f\|_{L^2(\R^d)}.
\end{split}
\]
The conclusion of Theorem~\ref{thm:RSI2} then follows from Theorem~\ref{thm:refined_dec}.
\end{proof}

To give some ideas about the proof of the refined decoupling inequality \eqref{eq:refined_dec}, note that heuristically $F 1_{Y_M} = \sum_{Q_j \subset Y_M} F 1_{Q_j}$ has Fourier support contained in an $R^{-1/2}$ neighbourhood of the paraboloid $\mathfrak{P}$. Thus we can apply the decoupling inequality \eqref{eq:BD} not at the finest scale, but at a coarser scale so that for some $K \leq R^{1/4}$,
\[
\|F\|_{L^p(Q_j)} \leq C_{\varepsilon/2}' K^{\varepsilon/2} \Big( \sum_{\beta \in P_{1/K}} \|F_{\beta}\|_{L^p(Q_j)}^2 \Big)^{1/2}
\]
for every $Q_j \subset Y_M$. Once we fix $Q_j$, for each $\theta \in P_{R^{-1/2}}$ there exists $O(1)$ many $T \in \W_{\theta}$ such that $Q_j \subset T$. Suppose for the moment, we have an extra assumption: 
suppose that for each $Q_j$, there exists $\sim M''$ $\beta$'s such that $F_{\beta} \ne 0$ on $Q_j$, and the number of $\theta$'s contained in each such $\beta$ so that $F_{\theta} \ne 0$ on $Q_j$ is $\sim M'$. 
Then by applying H\"{o}lder's inequality and taking $\ell^p$ norm over $Q_j \subset Y_M$, we have 
\[
\|F\|_{L^p(Y_M)} \leq C_{\varepsilon/2}' K^{\varepsilon/2}  (M'')^{\frac{1}{2} - \frac{1}{p}} \Big( \sum_{\beta \in P_{1/K}} \|F_{\beta}\|_{L^p(Y_M)}^p \Big)^{1/p},
\]
where 
\[
F_{\beta} := \sum_{\theta \in P_{R^{-1/2}}(\beta)} F_{\theta}.
\]
(Henceforth we write $P_{R^{-1/2}}(\beta)$ for a partition of $\beta$ into cubes of side length $R^{-1/2}$.)
Now we apply parabolic rescaling to $F_{\beta}$ that maps $\beta$ onto $[0,1]^d$. Wave packets at scale $R^{-1/2}$ that make up $F_{\beta}$ becomes wave packets at scale $R_1^{-1/2}$ that make up the rescaled $F_{\beta}$, where $R_1 := R/K^2$. A cube in $Y_M$ becomes part of a cube of side length $R_1^{1/2}$, which by our extra assumption has $M'$ wave packets from the rescaled $F_{\beta}$ passing through them. Thus $Y_M$ becomes a subset of $Y_{M'}$ for the rescaled $F_{\beta}$. Another way of saying the same thing is that we can first cover $\R^{d+1}$ by strips $S$ that are translates of $R^{1/2} K^{-1} \tau_{\beta}^*$. %
Here $\tau_{\beta}$ is defined as in \eqref{eq:tau_theta_def} (with $\beta$ in place of $\theta$) and $\tau_{\beta}^{\ast}$ is the dual parallelepiped defined in Section~\ref{sect:FT}.
The strips $R^{1/2} K^{-1} \tau_{\beta^*}$ have dimension $R^{1/2} \times \dots \times R^{1/2} \times R^{1/2} K$, and is thus the union of $K$ cubes of side length $R^{1/2}$. Since the tubes $T$ with $T \in \bigsqcup_{\theta \in P_{R^{-1/2}}(\beta)} \W_{\theta}$ are constrained in directions, the number of such $T$ passing through any of these $K$ cubes are comparable (in fact, $\sim M'$). Parabolic rescaling maps the strips $S$ into subcubes of $Q$ of side lengths $R_1^{1/2}$, and wave packets adapted to $T \in \bigsqcup_{\theta \in P_{R^{-1/2}}(\beta)} \W_{\theta}$ into wave packets adapted to tubes in $\T_{R_1^{-1/2}}$; see Figure~\ref{fig:Parabolic_rescaling} below.

\begin{figure}[ht] 
    \centering
    \begin{tikzpicture}
        \draw (0,0) to (3,0) to (4,5) to (1,5) to (0,0);
        \draw[blue] (0.8,0) to (1.6,0) to (3.6,5) to (2.8,5) to (0.8,0);
        \draw[thick] (2.2,3) to (2.8,3) to (3,4) to (2.4,4) to (2.2,3);
        \draw[thick] (0.8,2) to (1.4,2) to (1.6,3) to (1,3) to (0.8,2);
        \draw[ultra thick, ->] (4,2.5) to (6.5,2.5);
        \node at (2.6,3.5) {$S$};
        \node at (1.2,2.5) {$S$};
        \node[blue] at (2,2) {$T_{\theta,\nu}$};
        \node at (1.5,-0.4) {\large $\square_{\beta}$};
        \node at (5.25,2.9) {\footnotesize Parabolic rescaling};
        \draw (7,1.25) to (9.5,1.25) to (9.5,3.75) to (7,3.75) to (7,1.25);
        \draw[blue] (7+5/6*0.8,1.25) to (7+5/6*1.6,1.25) to (7+5/6*2.6,3.75) to (7+5/6*1.8,3.75) to (7+5/6*0.8,1.25);
        \draw[thick] (7+5/6*1.6,1.25+3/2) to (7+5/6*2.2,1.25+3/2) to (7+5/6*2.2,1.25+4/2) to (7+5/6*1.6,1.25+4/2) to (7+5/6*1.6,1.25+3/2);
        \draw[thick] (7+5/6*0.4,1.25+2/2) to (7+5/6*1,1.25+2/2) to (7+5/6*1,1.25+3/2) to (7+5/6*0.4,1.25+3/2) to (7+5/6*0.4,1.25+2/2);
        \node at (8.25,0.85) {\large $B_{R_1}$};
    \end{tikzpicture}
    \caption{In this figure we fix $\beta \in P_{1/K}$. On the left we have strips $S$ that are translates of $R^{1/2} K^{-1} \tau_{\beta}^*$. It is made up of $K$ cubes of side length $R^{1/2}$. Each wave packet from any $\theta \in P_{R^{-1/2}}(\beta)$ (pictured in blue) that passes through one of these $K$ cubes passes through all the others. Upon rescaling by a factor of $K$ in the short directions and $K^2$ in the long directions, the rescaled $S$ becomes a cube of side length $R_1^{1/2}$, and the rescaled wave packets becomes one at scale $R_1^{-1/2}$, with $R_1 := R/K^2$. If we consider the convex hull of all wave packets from $\theta \in P_{R^{-1/2}}(\beta)$ that passes through a given strip $S$, we get essentially a box $\square_{\beta}$, which gets rescaled into a cube $B_{R_1}$ of side length $R_1$.} 
    \label{fig:Parabolic_rescaling}
\end{figure}
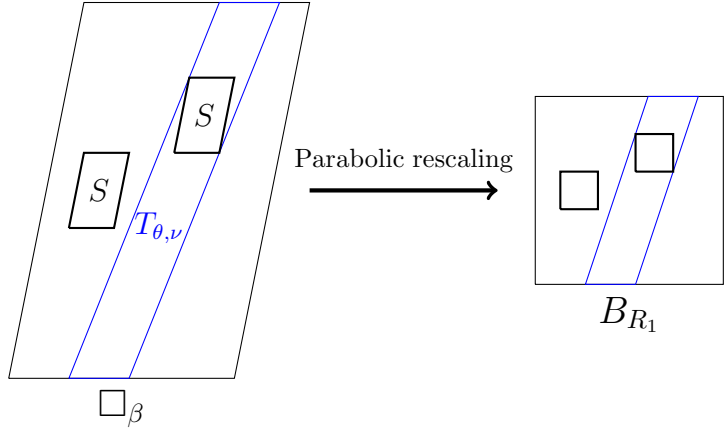
Thus applying our inductive hypothesis at the smaller scale $R_1$ as opposed to $R$, we have
\[
\|F_{\beta}\|_{L^p(Y_M)} \leq C_{\varepsilon} (R/K)^{\varepsilon} (M')^{\frac{1}{2}-\frac{1}{p}} \Big( \sum_{\theta \in P_{R^{-1/2}}(\beta)} \|F_{\theta}\|_{L^p(\R^{d+1})}^p \Big)^{1/p},
\]
which gives
\[
\|F\|_{L^p(Y_M)} \leq C_{\varepsilon/2}' K^{-\varepsilon/2}  C_{\varepsilon} R^{\varepsilon} (M'')^{\frac{1}{2} - \frac{1}{p}}   (M')^{\frac{1}{2}-\frac{1}{p}} \Big( \sum_{\theta \in P_{R^{-1/2}}} \|F_{\theta}\|_{L^p(\R^{d+1})}^p \Big)^{1/p}.
\]
Since %
$$
\#\{T \in \W : Q \subset T\} \sim M \iff \text{ there exists }\sim M \;\theta's \text{ such that }F_\theta \neq 0 \text{ on }Q.
$$
we have $M' M'' \leq M$. Thus the induction closes as long as $K = K(\varepsilon)$ is sufficiently large, so that $C_{\varepsilon/2}' K^{-\varepsilon/2} \leq 1$, and $C_{\varepsilon}$ was chosen large enough so that the induction hypothesis is true when say $R = 100$. All these were obtained under the above special assumption; but by pigeonholing, we can put ourselves under the special assumption by losing 2 powers of $\log R$. This gives a heuristic proof of Theorem~\ref{thm:refined_dec}; a more rigorous version is given below (it would be rigorous if we had perfect uncertainty principle, as in $\Q_p^{d+1}$).

\begin{proof}[Proof of Theorem~\ref{thm:refined_dec}]
We use a bootstrap argument. We know \eqref{eq:refined_dec} is true if $\varepsilon$ is replaced by a sufficiently large positive exponent. Suppose \eqref{eq:refined_dec} holds with $\eta > 0$ in place of $\varepsilon$. It suffices to show that \eqref{eq:refined_dec} also holds with $\frac{3}{4}\eta$ in place of $\varepsilon$, which we do as follows.

Let $K = R^{1/4}$. For $\beta \in P_{K^{-1}}$, let $\W_{\beta} := \bigsqcup_{\theta \subset \beta} \W_{\theta}$, and partition $\R^{d+1}$ into a disjoint union of strips $S$, that are translates of $R^{1/2} K^{-1} \tau_{\beta}^*$. For dyadic $M' \in [1,R^{d/2} K^{-d}]$ let $Y_{\beta,M'}$ be the union of all such strips $S$ for which $\#\{T \in \W_{\beta} \colon S \subset T\} \sim M'$. Since the strips have dimensions $R^{1/2} \times \dots \times R^{1/2} \times R^{1/2}K$, and $\R^{d+1} = \bigsqcup_{M'} Y_{\beta,M'}$ for every $\beta \in P_{K^{-1}}$, morally speaking, given $\beta \in P_{K^{-1}}$, each cube $Q$ with side length $R^{1/2}$ is contained in $Y_{\beta,M'}$ for a unique $M'$. In other words, for each $\beta \in P_{K^{-1}}$, we have
\[
\R^{d+1}= \bigsqcup_{M'} \{Q \colon Q \subset Y_{\beta,M'}\}.
\]
Thus we have
\[
F = \sum_{M'} \sum_{\beta \in P_{K^{-1}}} \sum_{Q \subset Y_{\beta,M'}} F_{\beta} 1_Q.
\]
Furthermore, let 
\[
N(Q,M') = \#\{\beta \in P_{K^{-1}} \colon Q \subset Y_{\beta,M'}\}
\]
and
\[
\mathcal{Q}_{M',M'',\beta} := \{Q \colon Q \subset Y_{\beta,M'}, N(Q,M') \sim M''\}.
\]
Then for each $\beta \in P_{K^{-1}}$, we have
\[
\R^{d+1}= \bigsqcup_{M',M''} \mathcal{Q}_{M',M'',\beta}.
\]
Hence
\[
F = \sum_{M',M''} \sum_{\beta \in P_{K^{-1}}} \sum_{Q \in \mathcal{Q}_{M',M'',\beta}} F_{\beta} 1_Q.
\]
For each cube $Q_0 \subset Y_M$ with side length $R^{1/2}$, we have
\begin{equation} \label{eq:dyadic_pigenhole}
    \|F\|_{L^p(Q_0)} \leq C (\log R)^2 \max_{M',M''} \Big\| \sum_{\beta \in P_{K^{-1}}} \sum_{Q \in \mathcal{Q}_{M',M'',\beta}} F_{\beta} 1_Q \Big\|_{L^p(Q_0)}
\end{equation}
so from the Bourgain-Demeter decoupling inequality \eqref{eq:BD}, that
\begin{equation} \label{eq:Q0_est}
\|F\|_{L^p(Q_0)} \leq C_{\varepsilon'} R^{\varepsilon'} (\log R)^2 \max_{M',M''} \Big( \sum_{\beta \in P_{K^{-1}}} \Big\| \sum_{Q \in \mathcal{Q}_{M',M'',\beta}} F_{\beta} 1_Q \Big\|_{L^p(Q_0)}^2 \Big)^{1/2}.
\end{equation}
Now fix dyadic $M', M''$. The $\ell^2$ norm above is zero unless $\Big\| \sum_{Q \in \mathcal{Q}_{M',M'',\beta}} F_{\beta} 1_Q \Big\|_{L^p(Q_0)} \ne 0$ for some  some $\beta \in P_{K^{-1}}$. In that case, we must have $Q_0 \in \mathcal{Q}_{M',M'',\beta}$, so $M' M'' \leq M$ by counting the number of $T \in \W$ for which $Q_0 \subset T$. This shows in \eqref{eq:Q0_est} we may restrict our attention to the dyadic numbers $M', M''$ for which $M' M'' \leq M$.
Also, using H\"{o}lder's inequality, we have
\[
\begin{split}
\Big( \sum_{\beta \in P_{K^{-1}}} \Big\| \sum_{Q \in \mathcal{Q}_{M',M'',\beta}} F_{\beta} 1_Q \Big\|_{L^p(Q_0)}^2 \Big)^{1/2} 
&= \Big( \sum_{\beta \in P_{K^{-1}}, Y_{\beta, M'} \supset Q_0} \Big\| \sum_{Q \in \mathcal{Q}_{M',M'',\beta}} F_{\beta} 1_Q \Big\|_{L^p(Q_0)}^2 \Big)^{1/2} \\
&\leq (M'')^{\frac{1}{2}-\frac{1}{p}} \Big( \sum_{\beta \in P_{K^{-1}}, Y_{\beta, M'} \supset Q_0} \Big\| \sum_{Q \in \mathcal{Q}_{M',M'',\beta}} F_{\beta} 1_Q \Big\|_{L^p(Q_0)}^p \Big)^{1/p} \\
&\leq  (M'')^{\frac{1}{2}-\frac{1}{p}} \Big( \sum_{\beta \in P_{K^{-1}}} \Big\| \sum_{Q \subset Y_{\beta,M'}} F_{\beta} 1_Q \Big\|_{L^p(Q_0)}^p \Big)^{1/p} 
\end{split}
\]
(The first inequality holds because $N(Q_0,M') \sim M''$, which follows from $Q_0 \in \mathcal{Q}_{M',M'',\beta}$ for some $\beta \in P_{K^{-1}}$.) Taking $\ell^p$ norm over all $Q_0 \subset Y_M$, this gives 
\[
\begin{split}
 \|F\|_{L^p(Y_M)} 
&\leq C_{\varepsilon'} R^{\varepsilon'} (\log R)^2 \max_{M',M'' \colon M' M'' \leq M} (M'')^{\frac{1}{2}-\frac{1}{p}} \Big( \sum_{\beta \in P_{K^{-1}}} \| F_{\beta} \|_{L^p(Y_{\beta,M'})}^p \Big)^{1/p} \\
&\leq  C_{\varepsilon'} R^{\varepsilon'} (\log R)^2 \max_{M',M'' \colon M' M'' \leq M} C_{\eta} (\frac{R}{K^2})^{\eta} (M'' M')^{\frac{1}{2}-\frac{1}{p}} \Big( \sum_{\beta \in P_{K^{-1}}} \sum_{\theta \in P_{R^{-1/2}}(\beta)} \|F_{\theta}\|_{L^p}^p \Big)^{1/p}
\end{split}
\]
where we applied a rescaled version of our bootstrap assumption. Hence we have
\[
\|F\|_{L^p(Y_M)} \leq C_{\varepsilon'} R^{\varepsilon'} (\log R)^2 C_{\eta} (\frac{R}{K^2})^{\eta} M^{\frac{1}{2}-\frac{1}{p}} \Big(\sum_{\theta \in P_{R^{-1/2}}} \|F_{\theta}\|_{L^p}^p \Big)^{1/p}.
\]
Plugging $K = R^{1/4}$ and taking $\varepsilon'$ sufficiently small shows that \eqref{eq:refined_dec} holds with $\frac{3}{4} \eta$ in place of $\varepsilon$ (in fact we can beat $R^{\eta}$ as long as we take $K$ to be any positive power of $R$ with exponent $\leq 1/4$). This closes the bootstrap argument.
\end{proof}

The above proof actually follows closely that in \cite{MR4055179}, except that we did not have to assume the non-zero $\|F_T\|_{L^p}$'s are comparable (which is possible because the right hand side of our inequality \eqref{eq:refined_dec} is an $\ell^p$ sum over $\theta$, whereas in \cite{MR4055179} the right hand side of their inequality involves an $\ell^2$ sum over the wave packets $T$), and that we did not need to pass through the $\square_{\beta}$'s pictured in Figure~\ref{fig:Parabolic_rescaling}; also we saved two pigeonholing steps that didn't seem essential. A similar proof can also be found in a recent preprint \cite{DORZ2023}, where they instead dyadically pigeonholed so that $\|F\|_{L^p(Q_0)}$ is comparable to a constant for a significant subset of cubes in $Y_M$. The subsequent pigeonholing in \cite{DORZ2023} relies on counting the number of these cubes to control $\|F\|_{L^p(Y_M)}$. On the contrary, we apply dyadic pigeonholing to the quantity on the right hand side of \eqref{eq:dyadic_pigenhole}, saving a pigeonholing step.

In practice one often needs a multilinear variant of the refined Strichartz inequality, or a fractal $L^2$ variant, as opposed to Theorem~\ref{thm:RSI} or \ref{thm:RSI2}; see e.g. \cite{MR3702674,MR3842310,hickman_survey}. One can either prove a multilinear refined Strichartz inequality by repeating the proof of its linear counterpart, as in \cite{MR3702674,MR3842310}, or by deducing it more directly from the linear refined decoupling estimate in Theorem~\ref{thm:refined_dec}, as in \cite{arXiv:2002.09525}. We omit the details.

\end{document}